\documentclass{article}
\usepackage{amsmath,amssymb,latexsym,amsfonts,amsthm} 
\usepackage[english]{babel}

\usepackage[normalem]{ulem} 

\usepackage{color}

\newtheorem{defn}{Definition}[section]
\newtheorem{lemma}[defn]{Lemma}
\newtheorem{ex}[defn]{Example}
\newtheorem{thm}[defn]{Theorem}
\newtheorem{prop}[defn]{Proposition}
\newtheorem{cor}[defn]{Corollary}
\newtheorem{rem}[defn]{Remark}

\newcommand{\h}{{H}}

\newcommand{\mn}{\mathbb N}
\newcommand{\mr}{\mathbb R}

\newcommand{\mc}{\mathbb C}

\newcommand{\me}{\mathcal{E}}
\newcommand{\mprop}{\mathcal{P}}
\newcommand{\sspace}{{\mathfrak{s}}}
\newcommand{\bs}{{\bf s}}

\newcommand{\spaceh}{H}
\newcommand{\nulel}{{\bf 0}}

\newcommand{\seqgr[1]}{(#1_n)_{n=1}^\infty }

\newcommand{\sech[1]}{{\cap_{k\in\mn_0}{#1}_k}}

\newcommand{\snorm[1]}{{\|\hspace{-.01in}|{#1}\|\hspace{-.015in}|}}

\def\be{\begin{equation}}
\def\ee{\end{equation}}
\def\newin {\kern-0.22em\in\kern-0.15em}
\def\newsubset {\kern-0.2em\subset\kern-0.2em}

\def\bp{\noindent{\bf Proof. \ }}
\def\ep{\noindent{$\Box$}}
\def\<{\langle}
\def\>{\rangle}

\title{Localization of Fr\'echet frames and expansion of generalized functions}
\author{Stevan Pilipovi\'c  and Diana T. Stoeva}
\date{\today}

\begin{document}
\maketitle

\begin{abstract}
 Matrix type operators with the off-diagonal decay of polynomial or sub-exponential 
  types are revisited with weaker assumptions concerning row or column estimates, still giving the continuity results for the frame type operators. Such results are extended from Banach to Fr\'{e}chet spaces. Moreover, the localization of  Fr\'{e}chet frames is used for the frame expansions of tempered distributions and a class of Beurling ultradistributions. 
\end{abstract}

\section{Introduction, Motivation and Main Aims}
Localized frames were introduced independently by
Gr\"{o}chenig \cite{G2004} and Balan, Casazza, Heil, and Landau \cite{BCHL1, BCHL2}.
The localization conditions in \cite{G2004} are related to 
off-diagonal decay (of polynomial or exponential type) of the matrix determined by the inner products of the frame elements and the elements of a given Riesz basis. A localized frame in this sense leads 
 to the same type of localization of the canonical dual frame as well as  to the convergence of the frame expansions in all  associated Banach spaces.
We refer to   
\cite{ChrStr, GC2004,FG, Futamura}, where various interesting properties and applications of localized frames were considered.
 The localization and self-localization, considered independently in \cite{BCHL1,BCHL2, BCHL}, are  directed to the  over-completeness of frames and the
relations between frame bounds and density with applications to Gabor frames.
For the present paper we have chosen to stick to the localization concept from \cite{G2004}, because the results obtained for a family of Banach spaces there can naturally be related to Fr\'echet frames (cf. \cite{pst,ps2,ps3,ps4}).

The aims of this paper are: 

First, to extend the continuity results on matrix-type operators acting on elements of a Banach or Fr\' echet spaces expanded by frames. 
Such results are related  to the off-diagonal decay conditions considered in the literature and the aim now is to use relaxed version of off-diagonal decay  which  requires column-decay but allows row-increase of a matrix. 

Second, to present the frame expansions of tempered distributions  $\mathcal S(\mathbb R)$ and tempered ultradistributions of Beurling type by the use of localization.

Beside of the main aims, the important novelty is the analysis related to sub-exponential off-diagonal decay without assumption  of the exponential off-diagonal decay as it was  considered in \cite{G2004}. 
More precisely, in \cite{G2004} the presumed exponential off-diagonal decay of matrices implies the analysis of  sub-exponentially weighted spaces.  Probably the most important impact in applications is related to Hermite basis
which is almost always used for the global expansion of $L^2$-functions or tempered generalized functions over ${\mathbb R}^n$. Our results by the use of localization, show that the same is true if one uses a kind of perturbation of Hermite functions through localization.

\label{key}
The paper is organized as follows.  
We recal in Section \ref{sdef}  the notation,   
basic definitions and the needed known results. 
In Section \ref{sec4} we consider matrices with column decay and possible row increase. For such type of matrices, we obtain  in Section \ref{sec_cont_frameoper}  continuity results for the frame related operators using less restrictive conditions in comparison to the localization conditions known in the literature. 
Sub-exponential localization is introduced and analyzed in Section \ref{sec3a}. The use of Jaffard's Theorem and \cite[Theorems 11 and 13]{G2004} is intrinsically connected with the sub-exponential localization. 
Section \ref{sec3} is devoted to Fr\'echet frames and series expansions in certain classes of Fr\'echet spaces based on polynomial, exponential, and sub-exponential localization. In particular, we obtain frame expansions in the Schwartz space $\mathcal{S}$ of rapidly decreasing functions and its dual, the space of tempered distributions, as well as in the spaces $\Sigma^\alpha$, $\alpha>1/2$, and their duals, spaces of tempered ultradistributions. 
In order to illustrate some results, we provide examples with the Hermite orthonormal bases $h_n, n\in\mn$, and construct a Riesz basis which is polynomially and exponentially localized to  $h_n, n\in\mn$. 
Finally, in the Appendix, we add some details in the proof of the  Jaffard's theorem.

\section{Notation, Definitions and Preliminaries}\label{sdef}
Throughout the paper, 
$(\h, \<\cdot,\cdot\>)$ denotes a separable Hilbert space and $G$ (resp. $E$) denotes the sequence $(g_n)_{n=1}^\infty$ (resp. $(e_n)_{n=1}^\infty$) with elements from $\h$.  
Recall that $G$ is called:

- \emph{frame for $\h$} \cite{DSframe}
 if there exist positive constants $A$ and $B$ (called {\it frame bounds}) so that 
$A\|f\|^2 \leq \sum_{n=1}^\infty |\<f,g_n\>|^2\leq B\|f\|^2$ for every $f\in\h$;

- \emph{Riesz basis for $\h$}  \cite{Bari51}
 if its elements are the images of the elements of an orthonormal basis under a bounded bijective operator on $\h$.

Recall  (see e.g. \cite{Cbook}), 
if $G$ is a frame for $\h$, 
then there exists a frame $\seqgr[f]$ for $\h$ so that 
$$f=\sum_{n=1}^\infty \<f,f_n\>g_n=\sum_{n=1}^\infty \<f,g_n\>f_n, f\in\h.$$ 
Such $\seqgr[f]$ is called a {\it dual frame of $\seqgr[g]$}. Furthermore, the {\it analysis operator} $U_G$, given by $U_Gf=(\<f, g_n\>)_{n=1}^\infty$, is bounded from $\h$ into $\ell^2$; the {\it synthesis operator} $T_G$, given by $T_Gf=\sum_{n=1}^\infty c_n g_n$, is bounded from $\ell^2$ into  $\h$;  
the {\it frame operator} 
$S_G:=T_G U_G$  is a bounded bijection of $\h$ onto $\h$ with unconditional convergence of the series $S_G f =\sum_{n=1}^\infty \<f, g_n\> g_n$. 
 The sequence $(S_G^{-1}g_n)_{n=1}^\infty$ is a dual frame of $\seqgr[g]$,  called the {\it canonical dual of $\seqgr[g]$}, and it will be denoted by 
 $(\widetilde{g_n})_{n=1}$ or  $\widetilde{G}$. 
When $G$ is a Riesz basis of $\h$ (and thus a frame for $\h$), then 
only $\widetilde{G}$ is a dual frame of $\h$, it is the unique biorthogonal sequence to $G$ and it is also a Riesz basis for $\h$. A frame $G$ which is not a Riesz basis has other dual frames in addition to the canonical dual and in that case we use notation $G^d$ or $(g_n^d)_{n=1}^\infty$ for a dual frame of $G$.

\vspace{.1in}
Next, $(X, \|\cdot\|)$ denotes a Banach space and 
 $( \Theta, \snorm[\cdot]) $ denotes a Banach sequence space; $\Theta$ 
 is called a {\it $BK$-space} if the coordinate functionals are continuous.
If the canonical vectors form a Schauder basis for $\Theta$, then $\Theta$ is called a {\it $CB$-space}. A $CB$-space is clearly a $BK$-space.  
 
Given a $BK$-space $\Theta$ and a frame $G$ for $\h$ with a dual frame $G^d=(g_n^d)_{n=1}^\infty$,  one associates to $\Theta$ the Banach space
$$\spaceh^{\Theta}_{G, G^d} := \{f\in\h \ : \ (\<f, g^d_n\>)_{n=1}^\infty \in\Theta, \, \|f\|_{\spaceh^{\Theta}_{G, G^d}}:=\snorm[ (\<f, g_n^d\>)_{n=1}^\infty]_\Theta\}.$$ 
When $G$ is a Riesz basis for $\h$, then 
we use notation $\spaceh^\Theta_G$ for  
 $\spaceh^{\Theta}_{G, \widetilde{G}}$.

 \subsection{Localization of frames}\label{localdef}
In this paper we consider polynomially and exponentially localized frames in the way defined in \cite{G2004}, and furthermore, sub-exponential localization. 
{\rm Let $G$ be a Riesz basis for the Hilbert space $\h$. A frame $E$ for $\h$ is called: 
\begin{itemize}
\item[- ] {\it polynomially localized with respect to $G$ with decay $\gamma>0$} 
 (in short, {\it $\gamma$-localized wrt $\seqgr[g]$})
if there is a constant $C_\gamma>0$ so that 
\begin{equation*}\label{poll}
\hspace{.5in}  \max\{ |\<e_m, g_n\>|, |\<e_m, \widetilde{g_n}\>|\} \leq 
 C_\gamma(1+|m-n|)^{-\gamma}, \ m,n\in\mn;
 \end{equation*} 
 \item[-] {\it exponentially localized with respect to $G$} if for some $\gamma>0$ there is a constant $C_\gamma>0$ so that
 \begin{equation*}\label{expl}
 \max\{ |\<e_m, g_n\>|, |\<e_m, \widetilde{g_n}\>|\} \leq 
 C_\gamma   \mathrm{e}^{-\gamma |m-n|}, \ m,n\in\mn.
 \end{equation*} 
\item[-] {\it $\beta$-sub-exponentially localized with respect to $G$} (for $\beta\in(0,1)$) if for some $\gamma>0$ there is  $C_\gamma>0$ so that
 \begin{equation*}
 \max\{ |\<e_m, g_n\>|, |\<e_m, \widetilde{g_n}\>|\} \leq 
 C_s   \mathrm{e}^{- \gamma |m-n|^\beta}, \ m,n\in\mn.
 \end{equation*}  
\end{itemize}
}

\subsection{Fr\'echet frames}\label{frechetfr}
We consider Fr\'echet spaces which are projective limits of Banach spaces as follows. 
Let $\{Y_k, | \cdot |_k\}_{k\in\mn_0}$ be a sequence of separable
Banach spaces such that 
\be \label{fx1} \{\nulel\} \neq
\sech[Y]\subseteq \ldots \subseteq Y_2 \subseteq Y_1 \subseteq Y_0, \; \ |\cdot|_0\leq | \cdot |_1\leq | \cdot |_2\leq
\ldots \ee 
\be \label{fx2} Y_F :=\sech[Y] \;\; \mbox{is dense in }
 Y_k, \;\; k\in\mn_0\;  (\mn_0=\mn\cup\{0\}). \ee

Under the conditions (\ref{fx1})-(\ref{fx2}), $Y_F$ is a Fr\'echet space.
 We will use such type of sequences in two cases:
1. $Y_k=X_k$ with norm $\|\cdot\|_k, k\in\mn_0;$
2. $Y_k=\Theta_k$ with norm $\snorm[\cdot]_k, k\in\mn_0$.

\vspace{.1in}
Let  $\{\Theta_k, \snorm[\cdot]_k\}_{k\in\mn_0}$ be a sequence of $CB$-spaces satisfying (\ref{fx1}). Then (\ref{fx2}) holds, because every sequence $(c_n)_{n=1}^\infty\in\Theta_F$ can be written as $\sum_{n=1}^\infty c_n \delta_n$ with the convergence in $\Theta_F$, where $\delta_n$ denotes the $n$-th canonical vector, $n\in\mn$.  
Furthermore,  $\Theta_F^*$ can be identified with the sequence space  $\Theta_F^{\circledast}:=\{(U\delta_n)_{n=1}^\infty \ : \ U\in \Theta_F^*\}$ with convergence naturally defined in correspondence with the convergence in $\Theta_F^*$.

\vspace{.1in}
 We use the therm \emph{operator} for a linear mapping, and by \emph{invertible operator on $X$} we mean a bounded bijective operator on $X$.  Given sequences of Banach spaces, $\{X_k\}_{k\in\mn_0}$ and $\{\Theta_k\}_{k\in\mn_0}$,
which satisfy (\ref{fx1})-(\ref{fx2}), an operator $T: \Theta_F \to X_F$
is called {\it $F$-bounded} if for every 
$k\newin\mn_0$, 
there exists a constant $C_k>0$ such that $\|T\seqgr[c]\|_k\leq
C_k \snorm[\{c_n\}_{n=1}^\infty]_k$ for all $\seqgr[c]\newin
\Theta_F$. 
\begin{defn}\label{fframe}  {\rm  \cite{ps4}} 
Let $\{X_k, \|\cdot\|_k\}_{k\in\mn_0}$ be a sequence of Banach
spaces satisfying (\ref{fx1})-(\ref{fx2}) and let
$\{\Theta_k, \snorm[\cdot]_k\}_{k\in\mn_0}$ be a sequence of
$BK$-spaces satisfying (\ref{fx1})-(\ref{fx2}). A sequence
$\seqgr[\phi]$ with elements from ${X_F^*}$ is called a {\it General  Fr\'echet frame} (in short, {\it General $F$-frame})  for $X_F$ with respect to $\Theta_F$ if 
there exist sequences $\{\widetilde{s}_k\}_{k\in\mn_0}\subseteq \mn_0$, $\{s_k\}_{k\in\mn_0}\subseteq \mn_0$, which increase to $\infty$ with the property $s_k\leq \widetilde{s}_k$, $k\in\mn_0$, and there
exist constants $0<A_k\leq B_k<\infty$, $k\in\mn_0$, satisfying 
\begin{equation}\label{fframestar}
(\phi_n(f))_{n=1}^\infty\in\Theta_F, \ f\newin X_F,
\end{equation}
\begin{equation}\label{fframetwostar}
A_k \|f\|_{s_k}\leq \snorm[\{\phi_n(f)\}_{n=1}^\infty]_k\leq
B_k\|f\|_{\widetilde{s}_k}, \ f\newin X_F, k\in\mn_0,
\end{equation}
and there exists a continuous operator 
$V:\Theta_F\rightarrow X_F$ so that $V(\phi_n(f))_{n=1}^\infty=f$
for every $f\newin X_F$.
 \end{defn}
 
 When $s_k=\widetilde{s}_k=k$, $k\in\mn_0$, and the continuity of $V$ is replaced by the stronger condition of $F$-boundedness of $V$, then the above definition reduces to the definition of a  {\it Fr\'echet frame} (in short, {\it $F$-frame})  for $X_F$ with respect to $\Theta_F$ introduced in \cite{ps2}. We will consider such frames in the sequel.
 
  In the particular case when $X_k=X$, and $\Theta_k=\Theta$,  $k\in\mn_0$,  a Fr\'echet frame  for $X_F$ with respect to $\Theta_F$ becomes a {\it Banach frame for $X$ with respect to $\Theta$}  as introduced in \cite{Gbanach}. 

For another approach to frames in Fr\'echet spaces we refer to \cite{BFGR}.  
For more on frames for Banach spaces, 
see e.g.  \cite{CHL, CCS, Spert} and the references therein.

\subsection{Sequence and function spaces} 
Recall that a positive
continuous 
function $\mu$ on $\mr$ is called: a  {\it $k$-moderate weight} if $k\geq 0$ and there exists a  constant $C>0$  so that
$ \mu(t+x)\leq C (1+|t|)^k \mu(x), \ t,x\in\mr;$
a  {\it $\beta$-sub-exponential} (resp. {\it exponential}) {weight}, if $\beta\in (0,1)$ (resp. $\beta=1$) and there exist  constants  $C>0,\gamma >0$, 
 so that
$ \mu(t+x)\leq C e^{\gamma |t|^\beta}\mu(x), \  t,x\in\mr.$ 
If $\beta$ is clear from the context, we will write just {\it sub-exponential weight}. 
Let 
 $\mu$ be a $k$-moderate, sub-exponential, or exponential weight so that $\mu(n)\geq 1$ for every $n\in\mn$, and $p\in[1,\infty)$. Then  
 the Banach space
$$\ell^p_{\mu}:=\{ (a_n)_{n=1}^\infty \ : 
\snorm[(a_n)_n]_{p,\mu} := (\sum_{n=1}^\infty |a_n|^p \mu(n)^{p})^{1/p}<\infty\}$$ 
is a  $CB$-space. 
We refer, for example, to  \cite[Ch. 27]{MV} for the  so called K\" othe sequence spaces. 
We will need the following, easy to prove, statements.
\begin{lemma} \label{lemdfr}
Let $G$ be a frame for $\h$ and let $G^d=\seqgr[g^d]$ be a dual frame of $G$. 
Let    $\mu_k$ be  $k$-moderate (resp. sub-exponential or exponential) weights, $k\in\mn_0,$ so that   
\begin{equation}\label{weightcond}
1=\mu_0(x)\leq \mu_1(x) \leq \mu_2(x) \leq ..., \ \mbox{ for every $x\in\mr$}. 
\end{equation}
Then the spaces $\Theta_k:=\ell^2_{\mu_k}$, $k\in\mn_0$, satisfy (\ref{fx1})-(\ref{fx2}). 
 Denote $M:= \{ (\<f,g^d_n\>)_{n=1}^\infty \ : \ f\in \h\}$. The  assumption that $M \cap \Theta_F$ 
is dense in  $M \cap \Theta_k\neq \{0\}$ with respect to the $\snorm[\cdot]_k$-norm  for every $k\in\mn$, 
leads to the conclusion that the spaces  $X_k:=\spaceh_{G,G^d}^{\Theta_k}$, $k\in\mn$, satisfy (\ref{fx1})-(\ref{fx2}).

If $G$ is a Riesz basis for $\h$, then   the density assumption of
 $M \cap \Theta_F$
 in  $M \cap \Theta_k\neq \{0\}$, $k\in\mn$, is fulfilled and in addition one has that $g_n\in X_F$ for every $n\in\mn$.
\end{lemma}
Throughout the paper we also consider specific weights, 
 relevant to the  function spaces of interest and the corresponding sequence spaces. 
Let  $k\in\mn_0$ and $\mu_k(x)=(1+|x|)^k$ (resp. $\mu_k(x)=e^{k|x|^\beta}, \beta\in(0,1]$), $ x\in\mr$.  Then, with  
$\Theta_k:=\ell^p_{\mu_k}$, $k\in\mn_0$, the projective limit $\cap_k \Theta_k$ 
is the space $\bs$ of rapidly  (resp. $\sspace^{\beta}$ of  sub-exponentially  when $\beta<1$ and exponentially when $\beta=1$) decreasing sequences determined by 
$ \{ \seqgr[a]\in\mc^\mn \ : \
 (\sum_{n=1}^\infty |a_n\mu_k(n))|^p)^{1/p} <\infty, \  \forall k\in\mn_0\}$,
which is the same set for any $p\in[1,\infty)$. The space $\bs$ (resp. $\sspace^\beta$) can also be derived as the projective limit of the 
Banach spaces
$\bs_k$ (resp. $\sspace^\beta_k$) defined as $\{(a_n)_{n=1}^\infty\in\mc^\mn  \ : 
\ \snorm[(a_n)_{n=1}^\infty]_{\sup,k}:= \sup_{n\in\mn} |a_n|n^k <\infty \}$, 
(resp. $ \snorm[\cdot]_{\sup,k}^\beta:= \sup_{n\in\mn} |a_n| e^{kn^{\beta}} <\infty \}$), $k\in\mn_0$;  
note that here instead of $k\in\{0,1,2,3,\ldots\}$ one can also use any strictly increasing sequence of non-negative numbers $k\in\{0, q_1, q_2, q_3, \ldots\}$.

 Recall that the well known Schwartz space {$\mathcal{S}$} is the intersection of Banach spaces
$${\mathcal S}_k (\mr) :=
\{f\in L^2(\mr):||f||_k=\sum_{m=0}^k||(1+|\cdot|^2)^{k/2}f^{(m)}||_{L^2(\mr)} \}, k\in\mn.$$
The dual $ {\mathcal S}' (\mr)$ is the space of tempered distributions.\\
\noindent
The space of sub-exponentially decreasing functions of order $1/\alpha$, $\alpha>1/2,$
is  $\Sigma^\alpha:=X_F=
\cap_{k\in\mn_0}\Sigma^{k,\alpha}$ where $\Sigma^{k,\alpha}$ are Banach spaces of $L^2-$ functions with finite norms
$$||f||_k^\alpha=\sup_{n\in\mn_0} ||\frac{k^{n}e^{k|x|^{1/\alpha}}|f^{(n)}(x)|}{ n!^\alpha}||_{L^2(\mr)}, k\in\mn.
$$
 Its dual $ ({\Sigma^\alpha}(\mr))'$ is the space of Beurling  tempered ultradistributions, cf.  \cite{sp,toft}.

\begin{rem}
		The case $\alpha=1/2$ leads to the trivial space $\Sigma^{1/2}=\{0\}.$ There is another way in considering the test space which corresponds to that limiting
	Beurling case $\alpha=1/2$ and can be considered also for $\alpha<1/2$ (cf. \cite{siam, toft,toft2,toft3}). 
	We will not treat these cases in the current paper.  
\end{rem}
In the sequel,  $(h_n)_{n=1}^\infty$ is
the Hermite orthonormal basis of $L^2(\mr)$ 
re-indexed from $1$ to $\infty$ instead of from $0$ to $\infty$. 
Recall that $h_n\in\Sigma^\alpha$, $\alpha>1/2$, $n\in\mn$.
Moreover, we know
\cite{Simon}:

- If $f\in\mathcal{S}$, then $(\<f, h_n\>)_{n=1}^\infty\in\bs$; conversely, if $(a_n)_{n=1}^\infty\in\bs$, then $\sum_{n=1}^\infty a_n h_n$ 
converges in $\mathcal{S}$ to some $f$ with $(\<f, h_n\>)_{n=1}^\infty=(a_n)_{n=1}^\infty$. 

-  If  $F\in \mathcal{S}'$, then 
$(b_n)_{n=1}^\infty:=(F(h_n))_{n=1}^\infty\in\bs'$ and $F(f)=\sum_{n=1}^\infty \<f, h_n\> b_n$,  $f\in  \mathcal{S}$; 
conversely, if $(b_n)_{n=1}^\infty\in\bs'$, then the mapping $F: f\to \sum_{n=1}^\infty \<f, h_n\> b_n$ 
 is well defined on $\mathcal{S}$, it determines $F$ as an element of $\mathcal{S}'$ 
and $(F(h_n))_{n=1}^\infty=(b_n)_{n=1}^\infty$.

The above two statements also hold when $\mathcal{S}$, $\mathcal{S}'$, $\bs$, and $\bs'$ are replaced by 
$\Sigma^\alpha$, $(\Sigma^\alpha)'$, $\sspace^{1/(2\alpha)}$, and $(\sspace^{1/(2\alpha)})'$ with $\alpha>1/2$, 
respectively (\cite{siam}, \cite{sp}, \cite{toft}).

We can consider  $\mathcal{S}$ and
$\Sigma^\alpha$ as the projective limit of Hilbert spaces $H^k, k\in\mn_0,$ with elements $f=\sum_n a_nh_n, $
in the first case with norms
$$\|f\|_{H^k}:=\snorm[(a_nn^k)_n]_{\ell^2}<\infty\}, k\in\mn_0,
$$ 
and in the second case with norms 
$$ \|f\|_{H^k}:=\snorm[(a_ne^{kn^{1/(2\alpha)}})_n]_{\ell^2}<\infty\}, k\in\mn_0.
$$
Thus, $(h_n)_n$ is an $F$-frame for $\mathcal S(\mr)$ with respect to $\bs$ as well as an $F$-frame for $\Sigma^\alpha$ with respect to $\sspace^{1/2\alpha}$, $\alpha>1/2$, 
($F$- boundedness is trivial).

\section{Matrix type operators}\label{sec4}

Papers \cite{G2004,GC2004,FG} concern matrices with off-diagonal decay of the form: for some $\gamma>0$ there is $C_\gamma>0$ such that 
\begin{equation}\label{ineqgrp} |A_{m,n}|\leq  \frac{C_\gamma}{(1+|m-n|)^\gamma} \ \ \  (\mbox{resp. } |A_{m,n}|\leq C_\gamma   \mathrm{e}^{-\gamma |m-n|}),\  \forall n,m\in\mn.
\end{equation}
In this paper we consider matrices with more general off-diagonal type of decay (see $ (***) $ below which is weaker condition compare to the polynomial type condition in (\ref{ineqgrp})). Moreover, we consider matrices which have column decrease but allow row increase   (see Propositions \ref{lem1} and \ref{propexp1alpha}) allowing sub-exponential type conditions as well.   
For such more general matrices, we generalize some results from \cite{G2004} with respect to certain Banach spaces and, furthermore, proceed to the Fr\'echet case. 

In the sequel, for a given matrix $(A_{mn})_{m,n\in\mn}$,  the letter $\mathcal{A}$ will denote the mapping $(c_n)_{n=1}^\infty \to (a_m)_{m=1}^\infty$ determined by $a_m=\sum_{n=1}^\infty A_{m,n} c_n$ (assuming convergence), $m\in\mn$; conversely, for a given mapping $\mathcal{A}$ 
determined on a sequence space containing the canonical vectors $\delta_n$, $n\in\mn$, 
the correspondent matrix $(A_{mn})_{m,n\in\mn}$ is given by $A_{m,n}=\langle \mathcal{A}\delta_n, \delta_m\rangle$. 
We will sometimes use $\mathcal{A}$ with the meaning of $(A_{mn})_{m,n\in\mn}$ and vice-verse.

\subsection{Polynomial type conditions}
Let us begin with some comparison of polynomial type of off-diagonal decay:
 
 \begin{lemma}\label{lem} Let $\gamma>0$. Consider the following conditions: 
$$ (*) \ \  \
 |A_{m,n}|\leq C \left\{ 
 \begin{array}{ll}
 \frac{(1+ m)^\gamma}{(1+n)^{2\gamma}}, & n\geq m,  \\
\frac{(1+ n)^\gamma}{(1+m)^{2\gamma}}, & n\leq m, \mbox{ for some } C>0.
\end{array}
 \right.
$$
$$(**) \ \ \  |A_{m,n}| \leq C(1+|n-m|)^{-\gamma},\mbox{ for some } C>0.$$
$$ (***) \ \ \ 
 |A_{m,n}|\leq C\left\{ 
 \begin{array}{ll}
 \frac{m^\gamma}{n^\gamma}, & n\geq m, \\
\frac{n^\gamma}{m^\gamma}, & n\leq m, \mbox{ for some } C>0.
\end{array}
 \right.
$$
Then, the implications $(*)\Rightarrow (**)\Rightarrow (***)$ hold. 
The converse implications are not valid. 
\end{lemma}
\bp Implications $(*)\Rightarrow (**)\Rightarrow (***)$ follow from the inequalities 
$\frac{(1+ min(m,n))^\gamma}{(1+max(m,n))^{2\gamma}}\leq 
 (1+|n-m|)^{-\gamma}
 \leq \frac{(min(m,n))^\gamma}{(max(m,n))^\gamma}, n,m\in\mn,$
which are easy to be verified. To show that 
$ (***)$ does not imply $(**)$ even up to a multiplication with a constant, take a matrix $A_{m,n}$ which satisfies 
$|A_{m,n}|=  \frac{C n^\gamma}{m^\gamma},\  n\leq m$, for some $\gamma>0$ and some positive constant $C,$
and assume that
there exist $\gamma_1(\gamma)\in\mn$ and a positive constant $K$ so that for $m\geq n$ one has
$\frac{C n^\gamma}{m^\gamma} \leq  \frac{K}{ (1+m-n)^{\gamma_1}}$; 
 then taking $m=2n$, one obtains $0<C \cdot 2^{-\gamma} \leq  \frac{K}{(1+n)^{\gamma_1}}\to 0 \mbox{ as $n\to\infty$,}$
which leads to a contradiction. 
In a similar spirit, one can show that 
$ (**)$ does not imply $(*)$. 
\ep

\vspace{.1in}
Below we show that the relaxed polynomial type conditions, as well as conditions allowing row-increase, still lead to continuous operators.

\begin{prop}\label{lem1}
Assume that the matrix
$(A_{mn})_{m,n\in\mn}$ satisfies the condition
$$
 |A_{m,n}|\leq \left\{ 
 \begin{array}{ll}
 C_0\, n^{\gamma_0}, & n>m, \\
  C_1 \,n^{\gamma_1} m^{-\gamma_1}, & n\leq m,
\end{array}
 \right.
$$
for some $\gamma_0\geq 0,  \gamma_1>0, C_0>0, C_1>0$.
Then  $\mathcal{A}$ is a continuous  operator 
from $\bs_{\gamma_1+\gamma_0+1+\varepsilon}$ 
 into $\bs_{\gamma_1}$ for any $\varepsilon\in(0,1)$.
\end{prop}
\bp 
Let $\varepsilon\in(0,1)$ and let $(c_n)_{n=1}^\infty\in \bs_{\gamma_0+\gamma_1+1+\varepsilon}$. 
For every  $n>m$, 
\begin{eqnarray*}
|A_{m,n} c_n| &\leq& C_0 |c_n| n^{\gamma_0}
\leq C_0  \left(\sup_j (|c_j| j^{\gamma_0+\gamma_1 +1+\varepsilon})\right) 
 \frac{1}{n^{\gamma_1+1+\varepsilon}}.
\end{eqnarray*} 
Next,
\begin{eqnarray*}
|\sum_{n=1}^m A_{m,n} c_n|
&\leq&  C_1 m^{-\gamma_1}
 \sum_{n=1}^m  |c_n| n^{\gamma_1}\\
&\leq&
 C_1 m^{-\gamma_1}   \snorm[(c_n)_{n=1}^\infty]_{\sup,\gamma_0+\gamma_1+1+\varepsilon} \sum_{n=1}^m \frac{1}{n^{\gamma_0+1+\varepsilon}}.
 \end{eqnarray*}
Therefore, 
\begin{eqnarray*}
|a_m| &\leq& |\sum_{n=1}^m A_{m,n} c_n|
  + |\sum_{n=m+1}^\infty A_{m,n} c_n| 
\\
&\leq& C_1 m^{-\gamma_1}   \snorm[(c_n)_{n=1}^\infty]_{\sup,\gamma_0+\gamma_1+1+\varepsilon} \sum_{n=1}^\infty \frac{1}{n^{\gamma_0+1+\varepsilon}} \\
&+& C_0   \snorm[(c_n)_{n=1}^\infty]_{\sup,\gamma_0+\gamma_1+1+\varepsilon}  
 \sum_{n=m+1}^\infty \frac{1}{ n^{\gamma_1+1+\varepsilon}}.
 \end{eqnarray*} 
Since  $\sum_{n=m+1}^\infty \frac{1}{ n^{\gamma_1+1+\varepsilon}}\leq m^{-\gamma_1} \sum_{n=m+1}^\infty \frac{1}{ n^{1+\varepsilon}},$ the
assertion follows.
\ep

\vspace{.1in}
A direct consequence of  Proposition \ref{lem1} is:
\begin{cor} \label{th1}
Assume that the matrix
$(A_{mn})_{m,n\in\mn}$ satisfies: there exist $\gamma_0\geq 0$ and $C_0>0$, and for every $\gamma>0$ there is $C_\gamma>0$ so that
$$
 |A_{m,n}|\leq \left\{ 
 \begin{array}{ll}
 C_0 n^{\gamma_0}, & n>m \\
  C_\gamma n^\gamma m^{-\gamma},  & n\leq m.
\end{array}
 \right.
$$
Then, $\mathcal{A}$  
is a continues operator from $\bs$ into $\bs$.
\end{cor}

In order to determine $\mathcal{A}$ as a mapping from a space $\bs_{\gamma_1}$ into the same space, we have to change the decay condition.  
\begin{prop} 
Let
$(A_{mn})_{m,n\in\mn}$ satisfy: 
$$(\exists \varepsilon>0, \gamma_1\in\mn) (\exists C_0>0, C_1>0) \mbox{ such that}$$
\begin{equation}\label{1newdecay}  
 |A_{m,n}|\leq \left\{ 
\begin{array}{ll}
C_0n^{-1-\varepsilon} , & n> m,  \\
C_1 n^{\gamma_1} m^{-\gamma_1-1-\varepsilon}, & n\leq m.
\end{array}
 \right.
\end{equation}
Then $\mathcal{A}$ is a continuous operator from
$\bs_{\gamma_1}$ into $\bs_{\gamma_1}$.
\end{prop}

\begin{rem}\label{Gcon}
For the same conclusion as above, one has in \cite{G2004} another condition non-comparible to (\ref{1newdecay}): 
\begin{equation}\label{G-decay}  
 |A_{m,n}|\leq C(1+|n-m|)^{-\gamma_1-1-\varepsilon}
 \end{equation}

\end{rem}

\subsection{Sub-exponential and exponential type conditions}\label{subsexp}

Up to the end of the paper $\beta$ will be a fixed number of the interval $(0,1];$  $\beta=1$ is related to the exponential growth order while
$\beta\in (0,1)$ corresponds to the pure sub-exponential growth order. 

\begin{prop} \label{propexp1alpha}  Assume that the matrix
$(A_{mn})_{m,n\in\mn}$ satisfies the condition: there exist positive constants $C_0, C_1$ and $\gamma_0\geq 0$, $\gamma_1>0$,  so that
\begin{equation}\label{cond1}
 |A_{m,n}|\leq \left\{ 
 \begin{array}{ll}
 C_0 \mathrm{e}^{\gamma_0n^{\beta}}, & n>m,  \\
C_1 \mathrm{e}^{-\gamma_1 (m^\beta-n^\beta)}, & n\leq m. 
\end{array}
 \right.
 \end{equation}
Then  $\mathcal{A}$ is a continues operator 
from $\sspace^\beta_{\gamma_1+\gamma_0+\varepsilon}$ 
into $\sspace^\beta_{\gamma_1}$ for any $\varepsilon\in(0,1)$.
\end{prop}
\bp  Let $\varepsilon\in(0,1)$ and let $(c_n)_{n=1}^\infty \in \sspace^\beta_{\gamma_1+\gamma_0+\varepsilon}$. Then for $n>m$,
\begin{eqnarray*}
|A_{m,n}c_{n}| &\leq&  
      C_0 |c_{n}| \mathrm{e}^{\gamma_0 n^\beta} 
\leq C_0  (\sup_j |c_j| \mathrm{e}^{(\gamma_1+\gamma_0+\varepsilon) j^\beta}) \mathrm{e}^{-(\gamma_1+\varepsilon) n^\beta}.
\end{eqnarray*} 
Further on,
\begin{eqnarray*}
|\sum_{n=1}^m A_{m,n} c_n| 
&\leq& C_1 \sum_{n=1}^m  |c_n| \mathrm{e}^{\gamma_1(n^\beta-m^\beta)}
\\
&\leq& C_1 \mathrm{e}^{-\gamma_1 m^\beta} 
(\sup_{j\in\mn} |c_j| \mathrm{e}^{(\gamma_1+\gamma_0+\varepsilon) j^\beta}) 
\sum_{n=1}^m \mathrm{e}^{-(\gamma_0 +\varepsilon)n^\beta}.
\end{eqnarray*}
Therefore, 
\begin{eqnarray*}
|a_m| 
&\leq& 
 \mathrm{e}^{-\gamma_1 m^\beta} 
\left(C_1 \sum_{n=1}^\infty \mathrm{e}^{-(\gamma_0+\varepsilon) n^\beta} 
+
 C_0 \sum_{n=1}^\infty  \mathrm{e}^{-\varepsilon n^\beta}\right)
 \snorm[(c_n)_n]_{\sup,\gamma_1+\gamma_0+\varepsilon}^\beta . 
\end{eqnarray*} \sloppy
This completes the proof.
\ep
 \begin{rem} Since $e^{-\gamma(m-n)^\beta} \leq e^{-\gamma(m^\beta-n^\beta)}$ for $n\leq m$ ($\beta\in(0,1], \gamma\in(0,\infty)$),  in (\ref{cond1}) we consider 
$e^{-\gamma(m^\beta-n^\beta)}$ instead of $e^{-\gamma{(m-n)^\beta}}.$
\end{rem}

\vspace{.1in} As a consequence of Proposition \ref{propexp1alpha}, we have:

\begin{cor} \label{propexp2alpha}  Assume that the matrix
$(A_{mn})_{m,n\in\mn}$ satisfies the condition: there exist  constants $C_0>0$ and $\gamma_0\geq 0$, and for every $\gamma>0$, there is a positive constant $C_\gamma$ so that
$$
 |A_{m,n}|\leq \left\{ 
 \begin{array}{ll}
 C_0 {e}^{\gamma_0n^{\beta}}, & n>m,  \\
C_\gamma e^{\gamma (n^\beta-m^\beta)}, & n\leq m. 
\end{array}
 \right.
 $$
Then $\mathcal{A}$ is a continuous operator from $\sspace^\beta$ into $\sspace^\beta$.
\end{cor}
\begin{prop}\label{subrem} 
Let
$(A_{mn})_{m,n\in\mn}$ satisfy the condition: there exist positive constants 
$\varepsilon, \gamma_1, C_0, C_1$, so that 
\begin{equation*}
 |A_{m,n}|\leq \left\{ 
\begin{array}{ll}
C_0e^{-\varepsilon n^\beta} , & n> m,  \\
C_1e^{\gamma_1 n^{\beta}} e^{-(\gamma_1 +\varepsilon)m^{\beta}}, & n\leq m.
\end{array}
 \right.
\end{equation*} Then $\mathcal{A}$ is a continues operator from
$\sspace^\beta_{\gamma_1}$ into $\sspace^\beta_{\gamma_1}$. 
\end{prop}

\begin{rem} 
One can simply show that the assumption  $|A_{m,n}|\leq C e^{-\gamma |m-n|^\beta}, m,n\in\mn$, leads to similar continuity results. We will consider this condition later
in relation to the  the invertiblity of such matrices and the Jaffard theorem. 
\end{rem}

\section{Continuity of the frame-related operators under relaxed ``decay'' conditions}
\label{sec_cont_frameoper}

We now determine weaker localization-conditions which  
are still sufficient to imply continuity of the frame-related operators.

\begin{prop} \label{pr1} 
Let $G$ be a frame for $\h$, $G^d$ be a dual frame of $G$, and
$ \mu_k(x) = (1+|x|)^k$, $k\in\mn_0$.
Under the notations in Lemma \ref{lemdfr}, assume that $M \cap \Theta_F$ 
is dense in  $M \cap \Theta_k\neq \{0\}$ with respect to the $\snorm[\cdot]_k$-norm  for every $k\in\mn$ and let  $E=(e_n)_{n=1}^\infty$ be a sequence with elements from $X_F$ which is a frame for $\h$. 
Then the following statements hold. 
  
 \begin{itemize}
\item[{\rm (i)}] 
 Assume that there exist $s_0\in\mn$, $C>0$
and for every $k\in\mn$ there exists  $C_k>0$  such that
\begin{equation*}
 |\<e_m, g_n\>| \leq 
\left\{ 
 \begin{array}{ll}
 Cn^{s_0}, & n>m, \\
 C_k n^k m^{-k}, & n\leq m.
\end{array}
 \right.
 \end{equation*} 
  Then the analysis operator $f\to U_Ef =(\<f,e_m\>)_{m=1}^\infty$ is continuous  one from  $X_F$ 
 into $\bs$.

\item[{\rm (ii)}] 
Assume that there exist $\widetilde{s}_0\in\mn_0$, $\widetilde{C}>0$ and for every $k\in\mn$ there exists  $\widetilde{C}_k>0$  
such that 
\begin{equation*}
  |\<e_m, g^d_n\>|\leq 
\left\{ 
 \begin{array}{ll}
 \widetilde{C}  m^{\widetilde{s}_0}, & m>n, \\
 \widetilde{C}_k m^k n^{-k}, & m\leq n.
 \end{array}
  \right.
 \end{equation*} 
  Then the synthesis operator $(c_n)_n\to T_E(c_n) =\sum c_n e_n$ is a continuous one from  $\bs$ 
 into $X_F$. 
 
\item[{\rm (iii)}] 
Under the assumptions of {\rm (i)} and {\rm (ii)},  
the frame operator $T_E U_E$ is continuous one from $X_F$ into  $X_F$. 
\end{itemize}
\end{prop}
\bp 
Note that under the given assumptions, $\Theta_F$ is the space $\bs$. 

(i) Let  $A_{m,n}=\< g_n, e_m\>$, $m,n\in\mn$, and  $\mathcal{A}$ be the corresponding operator for the matrix $A$. 
Let $f\in X_F$. Then $(\<f, g^d_n\>)_{n=1}^\infty\in\bs$ and 
$$\mathcal{A} (\<f, g^d_n\>)_{n=1}^\infty  = (\sum_{n=1}^\infty \< g_n, e_m\>\<f, g^d_n\>)_{m=1}^\infty=(\<f,e_m\>)_{m=1}^\infty.$$ By Corollary \ref{th1} it follows that $(\<f,e_m\>)_{n=1}^\infty\in\bs$. Furthermore, by Proposition \ref{lem1}, for every $k\in\mn$ there is a constant $K_{s_0,k,C,C_k}$ so that
\begin{eqnarray*}
\snorm[(\<f, e_m\> )_m]_{\sup,k} &=&\snorm[\mathcal{A}(\<f, g^d_n\> ))_n]_{\sup,k}\leq K_{s_0,k,C,C_k} \snorm[(\<f, g^d_n\> )_n]_{\sup,s_0 + k+2}\\
&\leq &K_{s_0,k,C,C_k} \, \snorm[(\<f, g^d_n\> )_n]_{\Theta_{k+2}} =K_{s_0,k,C,C_k} \, \|f\|_{s_0 + k+2}. 
\end{eqnarray*} 
Therefore, the analysis operator $U_E$ is continuous from  $X_F$  into $\bs$.

(ii)  Let $(c_n)\in \bs$. First we show that $\sum_{n=1}^\infty c_n e_n$ converges in $X_F$ and then the continuity of $T_E$. Since $(c_n)_{n=1}^\infty\in\ell^2$, $x=\sum_n c_n e_n\in\h$.  
Denote $A_{m,n}=\<  e_n, g^d_m\>$ and consider the corresponding operator $\mathcal{A}$. 
Then, $(\<x,g^d_m\>)_m=(\sum_n A_{m,n} c_n )_m = 
\mathcal{A}(c_n)\in\bs$, which implies that $x\in X_F$, and furthermore, for every $k\in\mn$, one has $\|T_E (c_n)_n\|_{k}= \|x\|_{k}=\snorm[(\<x, g^d_m\>)_m]_{\Theta_k}$. 
For every $k\in\mn$ there is a constant $R_k$ such that $\snorm[(d_n)]_{\Theta_k} \leq R_k \snorm[(d_n)]_{\sup,k+2}$ for every $(d_n)\in \bs_{k+2}$. By Proposition \ref{lem1}, we  conclude that
\begin{eqnarray*}
\|T_E (c_n)_n\|_{k}& \leq& 
R_k \snorm[(\<x, g^d_m\>)_m]_{\sup,k+2}=
R_k \snorm[\mathcal{A}(c_n)]_{\sup, k+2} \\
&\leq& R_k K_{(s_0, k, C, C_k)} \snorm[(c_n)]_{\sup,s_0 + k+4}. 
\end{eqnarray*}
 Thus, the synthesis operator $T_E$ is well defined and continues from  $\bs$  into $X_F$.

(iii) follows from (i) and (ii).
\ep

\vspace{.1in}

It is of interest to consider case when $X_F$ is $\mathcal{S}$. 
\begin{cor}
Let  $(e_n)_{n=1}^\infty$ be a frame of $L2(\mr)$  with elements in $\mathcal{S}(\mr)$. Assume that for every $k\in\mn$ there are constants $C_k, \widetilde{C}_k$ such that
\begin{equation}\label{egdecay7}
 |\<e_m, h_n\>| \leq 
\left\{ 
 \begin{array}{ll}
 C_k m^k n^{-k}, & n>m, \\
 \widetilde{C}_k n^k m^{-k}, & n\leq m.
\end{array}
 \right.
 \end{equation} 
 Then the analysis operator $U_E$ is continuous from $\mathcal{S}$ into $\bs$, the synthesis operator $T_E$ is  continuous  from $\bs$ into $\mathcal{S}$, and the frame operator $T_EU_E$ is  continuous  from $\mathcal{S}$ into $\mathcal{S}$. 
\end{cor}

Now, we consider sub-exponential 
weights.

\begin{prop} \label{pr2} 
Let $\beta\in(0,1)$ and let the assumptions of the first part of Lemma \ref{lemdfr} hold 
with the weights $ \mu_k(x) = e^{k|x|^\beta}$, $k\in\mn_0$.
Let 
  $E=(e_n)_{n=1}^\infty$ be a sequence with elements from $X_F$ which is a frame for $\h$.  
Then the following statements hold. 

 \begin{itemize}
\item[{\rm (i)}]
 Assume that there exist constants $\gamma_0\in\mn$, $C>0$ such that for every $k\in\mn$ there exists  $C_k>0$
such that 
\begin{equation}\label{egdecay5}
 |\<e_m, g_n\>| \leq 
\left\{ 
\begin{array}{ll}
 C \mathrm{e}^{\gamma_0n^{\beta}}, & n>m,  \\
C_k \mathrm{e}^{k (n^{\beta}-m^{\beta})}, & n\leq m, k\in\mn. 
\end{array}
 \right.
 \end{equation} 
Then the analysis operator $f\to U_Ef =(\<f,e_m\>)_{m=1}^\infty$ is a continuous one from  $X_F$ 
 into $\sspace^\beta$. 
 
\item[{\rm (ii)}]  
 Assume that there exist constants $\tilde \gamma_0\in\mn$, $\tilde C>0$ such that for every $k\in\mn$ there exists  $\tilde C_k>0$
such that
\begin{equation}\label{egdecay6}
 |\<e_m, g_n^d\>| \leq 
\left\{ 
\begin{array}{ll}
 \tilde C \mathrm{e}^{\tilde \gamma_0n^{\beta}}, & m>n,  \\
\tilde C_k \mathrm{e}^{k (m^\beta-n^\beta)}, & m\leq n. 
\end{array}
 \right.
 \end{equation} 
  Then the synthesis operator $(c_n)_n\to T_E(c_n) =\sum c_n e_n$ is a continuous one from  $\sspace^\beta$ 
 into $X_F$. 
 
\item[{\rm (iii)}] If (\ref{egdecay5}) and (\ref{egdecay6}) hold, then the frame operator $T_E U_E$ is  continuous  from $X_F$ into  $X_F$. 
\end{itemize}
\end{prop}
\bp Under the given assumptions, $\Theta_F$ is the space $\sspace^\beta$. The rest of the proof can be done in a similar way as the proof of Proposition \ref{pr1}, using Corollary \ref{propexp2alpha} instead of Corollary \ref{th1}. 
\ep

\vspace{.1in} 
If in the above proposition one chooses $G$ to be the Hermite basis $(h_n)_{n=1}^\infty$ and $\beta=1/(2\alpha)$, $\alpha>1/2$, then $X_F=\Sigma^\alpha$. 

\begin{cor} Let $\alpha>1/2$. 
Let $(e_n)_{n=1}^\infty$ be a sequence with elements from $\Sigma^\alpha$ which is a frame for $L^2(\mr)$ and such that for every $k\in\mn$ there are constants $C_k, \widetilde{C}_k$ such that
\begin{equation*}
 |\<e_m, h_n\>| \leq  C_k e^{-k|n^{1/(2\alpha)}-m^{1/(2\alpha)}|}, \ m,n,k\in\mn.
 \end{equation*} 
 Then the analysis operator $U_E$ is continuous from $\Sigma^\alpha$ into $\sspace^{1/(2\alpha)}$, the synthesis operator $T_E$ is  continuous  from $\sspace^{1/(2\alpha)}$ into $\Sigma^\alpha$, and the frame operator $T_EU_E$ is  continuous  from $\Sigma^\alpha$ into $\Sigma^\alpha$. 
\end{cor}

\section{Boundedness and Banach frames derived from sub-exponential localization of frames} \label{sec3a}

In this section we extend statements from \cite{G2004} for polynomially and exponentially localized frames to the case of sub-exponentially localized frames (Theorem \ref{thgr} below). 
We will use the Jaffard's theorem  \cite{sjaff} 
 given there for the sub-exponential and exponential case  (see Theorem \ref{sjt2} below).

First recall the {\it   Schur's test}: 
 {\it If $(A_{m,n})_{m,n\in\mn}$ is an infinite matrix satisfying \sloppy
$\sup_{s\in\mn} \sum_{t\in\mn} |A_{m,n}| \leq K_1$
and $\sup_{t\in\mn} \sum_{s\in\mn} |A_{m,n}| \leq K_2$, then the correspondent matrix frame type operator $\mathcal{A}$ 
is well defined and bounded from $\ell^p$ into $\ell^p$ for $1\leq p\leq \infty$ and the operator norm $\|A\|_{\ell^p\to\ell^p}\leq K_1^{1/p'} K_2^{1/p}$.  }

 Let $\beta\in (0,1]$ and $\gamma\in(0,\infty)$. Define $\me_{\gamma,\beta}$
 \label{eclass}
  to be the space 
 of matrices $(A_{m,n})_{m,n\in\mn}$ satisfying the following condition:
\begin{equation}\label{classdef}
 \exists C_\gamma\in (0,\infty) \ \mbox{so that} \ |A_{m,n}|\leq C_\gamma e^{-\gamma |m-n|^\beta}, m,n\in\mn.
\end{equation}

By the  Schur's test,  when $(A_{m,n})_{m,n\in\mn}\in\me_{\gamma,\beta}$, then 
the correspondent matrix type operator $\mathcal{A}$ 
is well defined and bounded from $\ell^2$ into $\ell^2$, and for the operator norm one has that $\|\mathcal{A}\|_{\ell^2\to\ell^2}\leq 2C_\gamma P_{\gamma,\beta}$, where $C_\gamma$ is the constant from (\ref{classdef}) and 
$P_{\gamma,\beta}$ denotes the sum of the convergent series $\sum_{j=0}^\infty e^{-\gamma j^{\beta}}$.

We will also need the following statements, which extend
 \cite[Lemmas 2 and 3]{G2004} to the case of sub-exponential localization.

\begin{lemma} \label{newsubexp} For every $\gamma\in(0,\infty)$ and $\beta\in (0,1]$, the following holds.

{\rm (a)} There exists a positive number $C$ so that 
$ \sum_{k\in\mn}  e^{-\gamma |m-k|^\beta} e^{-\gamma |k-n|^\beta} \leq C e^{-(\gamma/2) |m-n|^\beta}$ for every $s,t\in\mn$.

{\rm (b)} If the matrix $(A_{m,n})_{m,n\in\mn}$ belongs to $\me_{\gamma,\beta}$ and $\mu$ is a $\beta_\mu$-sub-exponential weight with $\beta_\mu<\beta$, then $\mathcal{A}$ maps boundedly $\ell^p_\mu(\mn_0)$ into $\ell^p_\mu(\mn_0)$   for every $p\in[1,\infty]$.

\end{lemma}

\bp
(a) can be proved following the idea of \cite[Lemma 2]{G2004}. 

(b) Let $\gamma_\mu$ comes from the $\beta_\mu$-sub-exponential weight $\mu$, i.e., 
$ \mu(t+x)\leq C_\mu e^{\gamma_\mu |t|^{\beta_\mu}}\mu(x), \  t,x\in\mr.$ 
 Take $k> \max(1, \gamma \beta /(\gamma_\mu \beta_\mu)) $ and use the assumption $\beta_\mu<\beta$ to observe that there is a constant $C_1\in(0,\infty)$ so that
$e^{-(\gamma/k) |x-n|^{\beta}}\mu(n)^{-1} \mu(x) \leq C_1 $. 
The rest of the proof can be done using a similar approach as in \cite[Lemma 3]{G2004}. 
\ep

 \begin{thm}{\rm (Jaffard)} \label{sjt2} 
 Given $\beta\in (0,1]$ and $\gamma\in(0,\infty)$,
let $(A_{m,n})_{m,n\in\mn}\in\me_{\gamma,\beta}$ and let the corresponding matrix type operator $\mathcal{A}$  be invertible on $\ell^2$. 
Then $\mathcal{A}^{-1}\in\me_{\gamma_1,\beta}$ for some $\gamma_1\in (0,\gamma)$. 
 \end{thm}

In the appendix we will give a sketch of the Jaffard's proof.

\begin{rem}
The exponential localization type condition $|A_{m,n}|\leq C_\gamma e^{-\gamma|m-n|}$, $m,n\in\mn$, considered in {\rm \cite{G2004}}, implies that 
$$\forall \beta\in(0,1) \ \forall k\in\mn \ \exists \widetilde{C_k} \mbox{  so that }|A_{m,n}|\leq \widetilde{C_k}e^{-k|s-t|^\beta}, m,n\in\mn.$$ 
Here we consider the more general case intrinsically related to  $\beta\in(0,1)$.
\end{rem}

\begin{thm} \label{thgr} 
Let $p\in[1,\infty)$ and $G$ be a Riesz basis for $\h$, and let $E$ be a frame for $\h$ which is  
$\beta$-sub-exponentially or exponentially localized (respectively, $(k+1+\varepsilon)$-localized for some $\varepsilon >0$) with respect to $G$. Let 
$\mu$ be $\beta_\mu$-sub-exponential weight 
and let $\beta_\mu <\beta$ in the case of $\beta$-sub-exponentially localized frame $E$
(respectively, let $\mu$ be a
$k$-moderate weight)  with $\mu(n)\geq 1$ for every $n\in\mn$.

Then for every $p\in[1,\infty)$ the following statements hold.
\begin{itemize}
\item[{\rm (i)}] 
 The analysis operator $U_E$ maps boundedly $\h_G^{\ell^p_\mu}$ into $\ell^p_\mu$.

\item[{\rm (ii)}]  The synthesis operator $T_E$ maps boundedly $\ell^p_\mu$ into $\h_G^{\ell^p_\mu}$.

\item[{\rm (iii)}] The frame operator $S_E=T_E U_E$ 
is invertible on $\h_G^{\ell^p_\mu}$ and the series in $S_Ef=\sum_{n=1}^\infty \<f, e_n\> e_n$ converges unconditionally.

\item[{\rm (iv)}] 
 The canonical dual frame $(\widetilde{e_n})_{n=1}^\infty$ of $(e_n)_{n=1}^\infty$ has the same type of localization as $(e_n)_{n=1}^\infty$, i.e., it is   
$\beta$-sub-exponentially or exponentially localized (resp. $(k+1+\varepsilon)$-localized) \mbox{
with respect to $G$.}

\item[{\rm (v)}] The frame expansions $f=\sum_{n=1}^\infty \<f,e_n\>\widetilde{e_n}=
\sum_{n=1}^\infty \<f,\widetilde{e_n}\>e_n$ hold with unconditional  convergence in $\h_G^{\ell^p_\mu}$.

\sloppy
\item[{\rm (vi)}] There is norm equivalence between 
$\|f\|_{\h_G^{\ell^p_\mu}}$,  $\snorm[(\<f,e_n\>)_{n=1}^\infty]_{\ell^p_\mu}$, and 
$\snorm[(\<f,\widetilde{e_n}\>)_{n=1}^\infty]_{\ell^p_\mu}$. 
\end{itemize}
\end{thm}

\bp In the cases of polynomial and exponential localization, the assertions are given in {\rm \cite[Prop. 8 and Prop. 10]{G2004}}.  For the sub-exponential case, one can proceed in a similar way, but to use the Jaffard's theorem Theorem \ref{sjt2} and Lemma \ref{newsubexp}. 
For the sake of completeness, we sketch a proof.

Consider the matrix $(A_{m,n})_{m,n\in\mn}$ with the property $|A_{m,n}|\leq Ce^{-\gamma|m-n|^\beta}$, for some $C>0$ and $\gamma>0$. 

(i) 
Let $f\in \h_G^{\ell^p_\mu}$
 and thus  $(|\<f, \widetilde{g_n}\>|)_{n=1}^\infty\in {\ell^p_\mu}$. 
 By Lemma \ref{newsubexp}(b), we have that $\mathcal{A} (|\<f, \widetilde{g_n}\>|)_{n=1}^\infty$
  belongs to ${\ell^p_\mu}$. Furthermore, for $m\in\mn$ we have
 \begin{eqnarray*}
 |\<f, e_m\>| 
 &\leq&
   \sum_{n=1}^\infty |\<f, \widetilde{g_n}\>\<g_n, e_m\>|
  \leq 
C \sum_{n=1}^\infty A_{m,n} |\<f, \widetilde{g_n}\>|. 
 \end{eqnarray*}
Therefore, $(\<f, e_m\>)_{m=1}^\infty$ also
 belongs to ${\ell^p_\mu}$ and 
 \begin{eqnarray*}  
\|(\<f, e_m\>)_{m=1}^\infty\|_{\ell^p_\mu} 
\leq C \|\mathcal{A} (|\<f, \widetilde{g_n}\>|)_{n=1}^\infty\|_{\ell^p_\mu} \leq 
 C\|\mathcal{A}\| \cdot \|f\|_{\h_G^{\ell^p_\mu}}.
 \end{eqnarray*}

(ii) Let $c=(c_n)_{n=1}^\infty\in \ell^p_\mu (\subseteq \ell^2)$. Then the series $\sum_{n=1}^\infty c_n e_n$ converges in $\h$ and let us denote its sum by $y$. 
Since $\mathcal{A}(|c_n|)_{n=1}^\infty\in \ell^p_\mu$ by Lemma \ref{newsubexp}(b),
and   since
$|\<y,  \widetilde{g}_m\>|\leq C \sum_{n=1}^\infty |c_n | |\<e_n, \widetilde{g}_m \> |
\leq C \sum_{n=1}^\infty A_{m,n} |c_n| $ for every $m\in\mn$, 
it follows that 
 $U_{\widetilde{G}} y \in \ell^p_\mu$ and therefore the element $T_G U_{\widetilde{G}} y=y$ belongs to $\h_G^{\ell^p_\mu}$. 
Hence, $T_E$ maps $\ell^p_\mu$ into $\h_G^{\ell^p_\mu}$
 and furthermore, 
 $\|T_E c\|_{\h_G^{\ell^p_\mu}}=\|y\|_{\h_G^{\ell^p_\mu}} 
 = \| U_{\widetilde{G}} y\|_{\ell^p_\mu}
 \leq C \| \mathcal{A} (|c_n|)\|_{\ell^p_\mu}$, which by Lemma \ref{newsubexp} implies that  
$\|T_E c\|_{\h_G^{\ell^p_\mu}}
 \leq C \| \mathcal{A}\|  \|c\|_{\ell^p_\mu}$.

(iii) 
By (i) and (ii), $S_E$ maps boundedly $\h_G^{\ell^p_\mu}$ into $\h_G^{\ell^p_\mu}$. For the unconditional convergence, 
take any re-ordering $N_1$ of $\mn$. 
Let $f\in\h_G^{\ell^p_\mu}$. Consider $\sum_{n\in N_1}\<f, e_n\> e_n$ and take $\varepsilon>0$. Since $(\<f, e_n\>)_{n\in N_1}\in \ell^p_\mu$,  there is a finite set $N_2\subset N_1$ so that  $\|(\<f, e_n\>)_{n\in N_1\setminus N_2}\|_{\ell^p_\mu} <\varepsilon$. 
Then for every finite $N_3$ such that $N_3\supset N_2$, $N_3\subset N_1$, one has that 
$\|Sf-\sum_{n\in N_3}\<f, e_n\> e_n\|
\leq   \|T_E\| \|(\<f, e_n\>)_{n\in N_1\setminus N_2}\|_{\ell^p_\mu} <  \|T_E\|\varepsilon.$ 
Therefore, $\sum_{n\in N_1}\<f, e_n\> e_n$ converges to $Sf$. 

Finally, let us show the bijectivity of $S_E$ on $\h_G^{\ell^p_\mu}$. Consider the operator $\mathcal{V}:= U_{\widetilde{G}} S_E T_G$ and observe that it is invertible on $\ell^2$ and it maps boundedly $\ell^p_\mu$ into $\ell^p_\mu$. 
Let $(V_{m,n})_{m,n\in\mn}$ be the corresponding matrix of $\mathcal{V}$. 
Since
$|V_{m,n}|  =|\<S_E g_n, \widetilde{g_m}\>|
\leq \sum_i |\< g_n, e_i\>|\cdot |\<e_i, \widetilde{g_m}\>|$, by Lemma \ref{newsubexp}(a), there is a positive constant $C$ so that
$|V_{m,n}|  \leq C e^{-(\gamma/2)|m-n|^\beta}$. 
Now by Theorem \ref{sjt2} it follows that $\mathcal{V}^{-1}\in \me_{\gamma_1,\beta}$  for some $\gamma_1\in (0, \gamma)$. 
By Lemma \ref{newsubexp}, it follows that $V^{-1}$ maps boundedly $\ell^p_\mu$ into $\ell^p_\mu$. 
Therefore, $V$ is a bounded bijection of $\ell^p_\mu$ onto $\ell^p_\mu$. 
  Now the representation $S_E=T_G V U_{\widetilde{G}}$ implies that $S_E$ is a bounded bijection of  $\h_G^{\ell^p_\mu}$ onto  $\h_G^{\ell^p_\mu}$.

(iv) For $m,n\in\mn$, 
 $\<\widetilde{e_m}, g_n\> = \sum_{j=1}^\infty \<e_m, g_j\> \overline{(V^{-1})_{nj}}$ and 
$\<\widetilde{e_m}, \widetilde{g_n}\> = \sum_{j=1}^\infty \<e_m, \widetilde{g_j}\> \overline{(V^{-1})_{nj}}$ and one can apply Theorem \ref{sjt2} and Lemma \ref{newsubexp}(a) to conclude.

(v) follows from (iii) 
and for (vi) one can use the representations 
$f= S_E^{-1} S_E =S_{\widetilde{E}^{-1}} S_{\widetilde{E}}f$ and 
the already proved (i)-(iv). 
\ep

\section{
Expansions in Fr\'echet spaces via localized frames} \label{sec3}

Our goal is expansion of elements of a Fr\'echet space and its dual via  localized frames and 
coefficients in a corresponding Fr\'echet sequence space. 
 First we present in the next theorem general results related to frames localized with respect to a Riesz basis. 
In the next secton we will  apply this theorem using frames localized with respect to the Hermite orthonormal basis in order to obtain frame expansions in the spaces $\mathcal{S}$ and $\Sigma^\alpha$, $\alpha>1/2$, and their duals. 
 
\begin{thm} \label{proprb1}
Let $G$ be a Riesz basis for $\h$,
 $k\in\mn_0$, 
and   $\mu_k$ be a $\beta_k$-sub-exponential (resp. $k$-moderate) weight so that    
(\ref{weightcond}) holds 
Let the spaces $\Theta_k$ and $X_k$ be as in Lemma \ref{lemdfr}. 
 Assume that $E=(e_n)_{n=1}^\infty$  is a sequence with elements in $X_F$ forming a frame for $\h$ which is $\beta$-sub-exponentially localized with $\beta>\beta_k$ for all $k\in\mn_0$ or exponentially  localized (respectively, $s$-localized $s$ for every $s\in\mn$) with respect to $G$.   
 Then, $\widetilde{e_n}\in X_F$, $n\in\mn$, and the following statements hold:
 \begin{itemize}
 
\item[{\rm (i)}] 
The analysis operator $U_E$ is $F$-bounded from $X_F$ into $\Theta_F$, the synthesis operator $T_E$ is $F$-bounded from $\Theta_F$ into $X_F$, and the frame operator $S_E$ is $F$-bounded and bijective from $X_F$ onto $X_F$ with unconditional convergence of the series in $S_Ef=\sum_{n=1}^\infty \<f, e_n\> e_n$. 
   \item[{\rm (ii)}] For every $f\in  X_F$, 
 \begin{equation}\label{reprxf}
 \mbox{$f=\sum_{n=1}^\infty \<f, \widetilde{e_n} \> e_n = \sum_{n=1}^\infty \<f, e_n\>\widetilde{e_n}$ (with convergence in $X_F$)}
 \end{equation}
with $(\<f,  \widetilde{e_n} \>)_{n=1}^\infty\in\Theta_F$ and $(\<f,  e_n \>)_{n=1}^\infty\in\Theta_F$.

\item[{\rm (iii)}] If $X_F$ and $\Theta_F$ have the following property with respect to $(g_n)_{n=1}^\infty$:
{\rm $$ \mbox{$\mprop_{(g_n)}$: For  $f\in \h$, one has $f\in  X_F$ if and only if $(\<f,  g_n \>)_{n=1}^\infty\in\Theta_F$.}$$}
 then $X_F$ and $\Theta_F$ also have the properties $\mprop_{(e_n)}$ and $\mprop_{(\widetilde{e_n})}$.

 \item[{\rm (iv)}]  
Both sequences  $({e_n})_{n=1}^\infty$ and $(\widetilde{e_n})_{n=1}^\infty$   form Fr\'echet frames for $X_F$ with respect to $\Theta_F$.

 \item[{\rm (v)}]  
  For every $g\in  X_F^*$, 
\begin{equation}\label{reprg}
\mbox{$g    = \sum_{n=1}^\infty g(e_n) \, \widetilde{{e_n}}
    =\sum_{n=1}^\infty g(\widetilde{e_n}) \, {e_n}$ (with convergence in $X_F^*$)} 
    \end{equation}
with $(g(e_n) )_{n=1}^\infty\in\Theta_F^*$ and $(g(\widetilde{e_n}) )_{n=1}^\infty\in\Theta_F^*$.

 \item[{\rm (vi)}]    If $(a_n)_{n=1}^\infty \in \Theta_F^*$, then $\sum_{n=1}^\infty a_n {e}_n$ (resp. $\sum_{n=1}^\infty a_n \widetilde{ e_n}$) converges in $X_F^*$, i.e., the mapping $f\mapsto \sum_{n=1}^\infty \<f, e_n\> a_n$  (resp. $f\mapsto \sum_{n=1}^\infty \<f, \widetilde{e_n}\> a_n$) determines a continuous linear functional on $X_F$.

\end{itemize} 

\end{thm}
\bp   
(i) The properties for $U_E$, $T_E,$ and $S_E$  follow easily using Theorem \ref{thgr}(i)-(iii). 

Further, the bijectivity of $S_E$ on $X_F$ implies that $\widetilde{e_n}\in X_F$ for every  $n\in\mn$.

(ii)  By 
 Theorem \ref{thgr}(v), for every $k\in\mn$ and every $f\in X_k$ we have that 
$f=\sum_{n=1}^\infty \<f, \widetilde{e_n} \> e_n = \sum_{n=1}^\infty \<f, e_n\>\widetilde{e_n}$  with convergence in $ X_k$. 
This implies that for every $f\in X_F$, one has that $ f=\sum_{n=1}^\infty \<f, \widetilde{e_n} \> e_n = \sum_{n=1}^\infty \<f, e_n\>\widetilde{e_n}$   with convergence in $X_F$.

 For every $k\in\mn$ and every $f\in X_k$, by 
 Theorem \ref{thgr}(i), we have that $(\<f, e_n\>)_{n=1}^\infty\in \Theta_k$. Therefore,  $(\<f, e_n\>)_{n=1}^\infty\in\Theta_F$ for every  $f\in X_F$.
  Furthermore, by 
  Theorem \ref{thgr}(iv), 
   $(\widetilde{e_n})_{n=1}^\infty$  
   has the same type of localization with respect to $G$ as $(e_n)_{n=1}^\infty$.
   Thus, applying 
   Theorem \ref{thgr}(i) 
    with $(\widetilde{e_n})_{n=1}^\infty$ as a starting frame, we get that $(\<f, \widetilde{e_n}\>)_{n=1}^\infty\in\Theta_F$ for $f\in X_F$.
 
(iii) 
If $f\in X_F$, it is already proved in (i) that $(\<f,  e_n \>)_{n=1}^\infty\in\Theta_F$ and $(\<f,  \widetilde{e_n} \>)_{n=1}^\infty\in\Theta_F$. 
To complete the proof of $\mprop_{(e_n)}$,  assume that $f\in\h$ is such that  $(\<f,  e_n \>)_{n=1}^\infty\in\Theta_F$. 
Consider 
$$
(\<f,  g_n \>)_{n=1}^\infty = (\<f, \sum_{j=1}^\infty \<g_n,\widetilde{e_j}\>e_j \>)_{n=1}^\infty 
=(\sum_{j=1}^\infty \<\widetilde{e_j},  g_n \>  \<f,e_j\>)_{n=1}^\infty. 
$$
Let $k\in\mn$. Since  $(\<f,  e_j \>)_{j=1}^\infty\in \ell^2_{\mu_k}$ and 
by Theorem \ref{thgr}(iv), $(\widetilde{e_n})_{n=1}^\infty$ has the same type of localization with respect to $G$ as $(e_n)_{n=1}^\infty$, it follows from 
 Lemma \ref{newsubexp}(b)  (for the case of sub-exponential localization) and from the way of the proof of \cite[Lemma 3]{G2004} (for the case of polynomial and exponential localization) 
that $(\sum_{j=1}^\infty \<\widetilde{e_j},  g_n \>  \<f,e_j\>)_{n=1}^\infty \in \ell^2_{\mu_k}$. 
Therefore, $(\<f,  g_n \>)_{n=1}^\infty\in\Theta_F$ and thus, by $\mprop_{(g_n)}$, it follows that $f\in X_F$. 
For completing the proof of $\mprop_{(\widetilde{e_n})}$, if $f\in \h$ is such that $(\<f,  \widetilde{e_n} \>)_{n=1}^\infty\in\Theta_F$, it follows in a similar way as above that $f\in X_F$.

(iv) 
 By (i),  
  $(\bold{e}_n(f))_{n=1}^\infty\in\Theta_F$
  for $f\in X_F$, and by 
Theorem \ref{thgr}(vi), 
for $k\in\mn$ and  $f\in X_k$, the norms $\snorm[(\<f,e_n\>)_{n=1}^\infty]_{\Theta_k}$ and $\|f\|_{X_k}$ 
 are equivalent. 
 Furthermore,  it follows from Theorem \ref{thgr} that the operator $V:=S_E^{-1} T_E\vert_{_{\Theta_F}} $ maps $\Theta_F$ into $X_F$ and it is $F$-bounded. 
Clearly,  $V(\bold{e}_n(f))_{n=1}^\infty=f$, $f\in X_F$. 
 Therefore, $(\bold{e}_n)_{n=1}^\infty$ is an $F$-frame for $X_F$ with respect to $\Theta_F$. 
 In an analogue way,  $(\widetilde{\bold e_n})_{n=1}^\infty$ is also an $F$-frame for $X_F$ with respect to $\Theta_F$.

(v)  The representations in (i) can be re-written as 
$f=\sum_{n=1}^\infty \widetilde{{\bf e}_n} (f) e_n = \sum_{n=1}^\infty  {\bf e}_n (f) \widetilde{e_n}$, $f\in X_F$, which implies validity of (\ref{reprg}) for $g\in X_F^*$. 

For the rest of the proof, consider the $F$-bounded (and hence continuous) operator $V$ from the proof of (iv) and observe that $\widetilde{e_n}=V\delta_n$, $n\in\mn$. This implies that for $g\in X_F^*$ we have
$(g(\widetilde{e_n}))_{n=1}^\infty =(gV (\delta_n))_{n=1}^\infty\in \Theta_F^\circledast$. With similar arguments, considering the operator $\widetilde{V}=S_{\widetilde{E}}^{-1} T_{\widetilde{E}}\vert_{_{\Theta_F}} $, it follows that 
$(g(e_n))_{n=1}^\infty \in \Theta_F^\circledast$.

(vi) Let $(a_n)_{n=1}^\infty \in \Theta_F^*$ and thus there is $k_0\in\mn$ so that $(a_n)_{n=1}^\infty \in \Theta_{k_0}^*$, i.e., $C:=\sum_{n=1}^\infty |a_n|^2 |\mu(n)|^{-2k_0}<\infty$. 
By Theorem \ref{thgr}(vi), there is a positive constant $B_{k_0}$ so that 
$\snorm[(\<f, e_n\>)_{n=1}^\infty]_{\Theta_{k_0}} \leq B_{k_0} \|f\|_{X_{k_0}}$ for every $f\in X_F$.  
Let $f\in X_F$. By (i), 
$(\<f,  e_n \>)_{n=1}^\infty\in \Theta_F$. 
Therefore,  $\sum_{n=1}^\infty \<f, e_n\> a_n$ converges 
and furthermore, 
\begin{eqnarray*}
|\sum_{n=1}^\infty \<f, e_n\> a_n|^2
&\leq&
 \left(\sum_{n=1}^\infty |\<f, e_n\>|^2 |\mu(n)|^{2k_0}\right)
 \left(\sum_{n=1}^\infty |a_n|^2 |\mu(n)|^{-2k_0}\right) \\
 &=&
 C \snorm[(\<f, e_n\>)_{n=1}^\infty]_{\Theta_{k_0}}
 \leq B_{k_0} C   \|f\|_{X_{k_0}}, 
 \end{eqnarray*}
which implies continuity of the linear mapping $f\mapsto \sum_{n=1}^\infty \<f, e_n\> a_n$. 
In a similar way, it follows that  $f\mapsto \sum_{n=1}^\infty \<f, \widetilde{e_n}\> a_n$ determines a continuous linear functional on $X_F$. 
\ep

\vspace{.1in}

\begin{rem} Note that in the setting of the above theorem, when $G$ is an orthonormal basis of $\h$ or more generally, when $G$ is a Riesz basis for $\h$  satisfying any of the following two conditions:

 $(\mathcal{P}_1)$:\ \ $\forall s\in\mn \ \,  \exists C_s>0 \ : \  |\<g_m, g_n\>|\leq  C_s(1+|m-n|)^{-s}, \ m,n\in\mn$,

 $(\mathcal{P}_2)$:\ \ $\exists s>0 \ \, \exists C_s>0 \ : \    |\<g_m, g_n\>| \leq 
 C_s   \mathrm{e}^{-s |m-n|}, \ m,n\in\mn,
$

\noindent
 then the property $\mprop_{(g_n)}$ is  satisfied. 

\end{rem}

\subsection{Frame expansions of tempered distributions and ultradistributions }

Here we apply Theorem \ref{proprb1} to obtain series expansions in the spaces $\mathcal{S}$ and $\Sigma^\alpha$ ($\alpha>1/2$), and their dual spaces, via frames which are localized with respect  to the Hermite basis.

\begin{thm} \label{prophb}  
 Assume that the sequence $(e_n)_{n=1}^\infty$  with elements from $\mathcal{S}(\mr)$ (resp. in $\Sigma^\alpha$) is a  frame for $L^2(\mr)$ which is polynomially  (resp. sub-exponentially or exponentially) localized with respect to the Hermite basis $(h_n)_{n=1}^\infty$ with decay $\gamma$ for every $\gamma\in\mn$. Let $(g_n)_{n=1}^\infty=(h_n)_{n=1}^\infty$. 
  Then $\mprop_{(g_n)}$ and the conclusions in Theorem \ref{proprb1} hold with $X_F$ replaced by $\mathcal{S}$ 
 (resp. $\Sigma^\alpha$) and $\Theta_F$ replaced by $\bs$ (resp.  $\sspace^{1/(2\alpha)}$). 
  \end{thm}
\bp
 For $k\in\mn_0$, consider the $k$-moderate weight $ \mu_k(x) = (1+|x|)^k$
 the spaces $\Theta_k:=\ell^2_{\mu_k}$, $k\in\mn_0$, satisfy (\ref{fx1})-(\ref{fx2}) and their projective limit 
 $\Theta_F$ is the space $\bs$. 
Consider the spaces   $X_k:=\spaceh^{\Theta_k}_{(h_n)}$, $k\in\mn_0$, which satisfy (\ref{fx1})-(\ref{fx2}). 
As observed after Theorem \ref{proprb1}, the property $\mprop_{(h_n)}$ is satisfied. 
Since for  $f\in L^2(\mr)$ one has that $f\in\mathcal{S}$ if and only if $(\<f,h_n\>_{n=0}^\infty)\in\bs$, it now follows that $X_F=\mathcal{S}$. Then the conclusions of Theorem \ref{prophb} follow from Theorem \ref{proprb1}. 

The respective part of the theorem follows in a similar way. 
\ep

\vspace{.1in}
As noticed in \cite{ps5}, 
having in mind the known  expansions of tempered distributions $(\mathcal S(\mr_+))'$ \cite{GTeissier,Duran} 
 and 
Beurling ultradistributions $(G^\alpha_\alpha(\mr_+))'$ \cite{Duran2,JPP,bpsj}, 
and their test spaces, by the use of the Laguerre orthonormal basis $l_n, n\in\mn,$ and validity of the corresponding properties $\mprop_{(l_n)}$, 
we can transfer the above results to the mentioned classes of distributions and ultradistributions over $\mr_+.$

\begin{rem} \label{rempol} For Proposition \ref{pr1} (resp. \ref{pr2}), it is of interest to consider cases when $X_F$ is the space  $\mathcal{S}$ (resp. $\Sigma^\alpha$). Based on Theorem \ref{prophb}, we can clarify such cases. 
If  $G$ is a frame for $L^2(\mr)$ with elements from $\mathcal{S}$ and  polynomially localized with respect to $(h_n)_{n=1}^\infty$, then one also has $X_F=\mathcal{S}$  since in this case Theorem \ref{prophb} 
 implies that $f\in \mathcal{S}$ if and only if $(\<f, \widetilde{g_n}\>)\in \bs$, and besides that one also has that $f\in X_F$ if and only if 
$(\<f, \widetilde{g_n}\>)\in \bs$. 

Concerning Proposition \ref{pr2}, with $\alpha>1/2$ and $\beta=1/(2\alpha)$, we have the similar conclusion for $\Sigma^\alpha$.
\end{rem}

\begin{ex}  \label{expol} As an illustration of Theorems \ref{proprb1} and \ref{prophb}, we give the next example.
 Let $r\in\mn$ and for $i=1,2,\ldots,r$, take  $\varepsilon_i \geq 0 $  and a sequence $(a_n^i)_{n=1}^\infty$ of complex numbers
 satisfying $|a_n^i|\leq \varepsilon_i$ for $n\geq 2$,  $\sum_{i=1}^r |a^i_1|\leq 1$, and $\sum_{i=1}^r \varepsilon_i <1$. For $n\in\mn$, consider $e_n:=h_n + \sum_{i=1}^r a_n^i h_{n+i}$. The sequence $(e_n)_{n=1}^\infty$ is a Riesz basis for 
$L^2(\mr)$ and it is $s$-localized with respect to the Hermite orthonormal basis $(h_n)_{n=1}^\infty$ for every $s>0$, as well as exponentially localized with respect to $(h_n)_{n=1}^\infty$. 
In order to show that $(e_n)_{n=1}^\infty$ is a Riesz basis for $L^2(\mr)$, we will represent $(e_n)_{n=1}^\infty$ as a sequence $(Uh_n)_{n=1}^\infty$ for some bounded bijective operator from $L^2(\mr)$ onto $L^2(\mr)$
using similar techniques as in \cite[Example 1]{CC}.  
 Define $U$ by $Uh_n:=h_n + \sum_{i=1}^r a_n^i h_{n+i}$, $n\in\mn$, and by linearirty on the linear span of $h_n, n\in\mn$. 
The obtained operator is bounded on the linear span, so extend it by continuity on $L^2(\mr)$, 
leading to a bounded operator on $L^2(\mr)$. 
It remains to show the bijectivity of $U$.  
Let  $f\in L^2(\mr)$ Then 
$$ Uf - \<f,h_1\>h_1  =\sum_{n=2}^\infty \<f,h_n\> h_n 
+ \sum_{n=1}^\infty \sum_{i=1}^r a^i_n  \<f,h_n\>h_{n+i} 
\in \overline{span}\{h_n\}_{n=2}^\infty
$$  and thus,
$\|Uf \|\geq |\<f,h_1\>|$. 
Furthermore, for $f\in L^2(\mr)$, 
     $$Uf=\sum_{n=1}^\infty \<f,h_n\> Uh_n 
   =f+ \sum_{n=1}^\infty \sum_{i=1}^r a^i_n  \<f,h_n\>h_{n+i}, $$
   leading to        
              \begin{eqnarray*}\|Uf\|&\leq& \|f\| + \sum_{i=1}^r
\|\sum_{n=1}^\infty a^i_n \<f,h_n\>h_{n+i}\|
=\|f\|+\sum_{i=1}^r
\|(a^i_n \<f,h_n\>)_{n=1}^\infty\|_{\ell^2}\\
&\leq&
 \|f\|+\sum_{i=1}^r \sqrt{|a^i_1|^2 \|f\|^2 + \varepsilon_i^2 \|f\|^2}
\leq 
\|f\|+\sum_{i=1}^r (|a^i_1|  + \varepsilon_i) \|f\|)\\
&\leq& 3\|f\|,
\end{eqnarray*}

\begin{eqnarray*}\|Uf-f\|&\leq& 
\sum_{i=1}^r
\sqrt{|a^i_1|^2 |\<f,h_1\>|^2 + \sum_{n=2}^\infty|a^i_n|^2 |\<f,h_n\>|^2 } \\
 &\leq& \|Uf\| +  (\sum_{i=1}^r  \varepsilon_i) \|f\|\\
&\leq& 
 \frac{3+\sum_{i=1}^r  \varepsilon_i}{4}\|Uf\| +\frac{3+\sum_{i=1}^r  \varepsilon_i}{4}\|f\|.
\end{eqnarray*}
Since $ \frac{3+\sum_{i=1}^r  \varepsilon_i}{4}<1$, it follows from \cite[Lemma 1]{CC} that the bounded operator $U$ is bijective on $L^2(\mr)$ and thus $(e_n)_{n=1}^\infty$ is a Riesz basis for $L^2(\mr)$.

Note that under the assumptions of the example, the classical way to obtain invertibility 
of $U$ does not apply, because $\|U-Id\|$ is not necessarily  smaller then $1$. 
 \end{ex}

\begin{rem} As explained in \cite{G2004}, when dealing with Gabor frames and localization, the natural bases to be considered in this respect are the Wilson bases, but then Gabor frames are not localized with respect to a Wilson basis in the strict sense of the definition of polynomial and exponential localization. However, under appropriate conditions, the authors of \cite{G2004} still obtain statements in the spirit of 
 Theorem \ref{thgr} 
which now leads to conclusions as in Theorem \ref{proprb1}.
\end{rem}

\section{Appendix} 

We add in Jaffard's proof some comments in the end. 
We also give a simple known assertions 
for the class $\me_{\gamma,\beta}$ defined on page \pageref{eclass}. 
   \begin{lemma}\label{ab} Let $\mathcal{A}\in\me_{\gamma,\beta}$. 
   \begin{itemize}
   \item[{\rm (i)}]  If $B$ belongs to $\me_{\gamma',\beta}$ for some $\gamma'\in(0,\gamma)$, then $AB$ belongs to $\me_{\gamma',\beta}$. 
   \item[{\rm (ii)}] If $B$ belongs to $\me_{\gamma,\beta}$, then $AB$ belongs to $\me_{\gamma',\beta}$ for every $\gamma'\in(0,\gamma)$. 
   \end{itemize}
   \end{lemma}
 
 \bp 
(i) Under the assumptions, there exist positive constants $C_A$ and $C_B$ so that 
 $\ |A_{m,n}|\leq C_A e^{-\gamma |m-n|^\beta}$ and 
  $\ |B_{m,n}|\leq C_B e^{-\gamma' |m-n|^\beta}$ for $m,n\in\mn$. Further, using $|m-k|^{\beta}+|k-n|^\beta\geq |m-n|^\beta,$ for every $m,n\in\mn$, one has
$$
|(AB)_{m,n}|
\leq C_A C_B \sum_{k=1}^{\infty} 
 e^{-\gamma'|m-k|^\beta}e^{-\gamma|k-n|^\beta}
\leq
 C_A C_B   e^{-\gamma' |m-n|^\beta}  2\sum_{j=0}^{\infty}
 e^{-(\gamma-\gamma') j^\beta}.
 $$ 
(ii) follows from (i).
 \ep

 \vspace{.1in}
 Now we give some details for the Jaffard's proof of Theorem \ref{sjt2}, providing explicit estimates for the bounds.
 
    As in \cite{sjaff}, $AA^*$  belongs to $\me_{\gamma',\beta}$ for $\gamma'\in(0,\gamma)$  by Lemma \ref{ab}, $AA^*=\|AA^*\|(Id -R)$ for some operator $R$ with $ ||R||=r<1$, and  $(AA^*)^{-1}=\|AA^*\|^{-1} \sum_{n=0}^\infty R^n$. 
    With the method from \cite{sjaff}, 
 one obtains     
$$ \sum_{k=1}^\infty  |(R^k)_{m,n}|\leq 
 \sum_{k=1}^\infty \min\{r^k,  (C_1)^k (2P)^{k-1}e^{\gamma'' |m-n|^\beta}\}
$$
where $\gamma''\in(0,\gamma')$, $C_1=1+\frac{C_{AA^*}}{\|AA^*\|}$ (with $C_{AA^*}$ satisfying $\ |(AA^*)_{m,n}|\leq C_{AA^*} e^{-\gamma' |m-n|^\beta}$ for every $s,t\in\mn$),     and $P=P_{\gamma'-\gamma'',\beta}=
\sum_{j=0}^\infty e^{-(\gamma'-\gamma'')j^{\beta}}.$
 
 For k=0,  $(R^k)_{m,n}$ equals $1$ when $m=n$ and $0$ otherwise, so we have 
 $ |(R^0)_{m,n}|\leq  e^{-p |m-n|^\beta}$ for every $p$. 
 
Fix $|m-n|$ and let $\varepsilon\in(0,1)$. With $K:=C_1 2P (>1)$, let 
 \begin{eqnarray*}
 \Lambda_{|m-n|}&:=&\{n\in\mathbb N \  : \ K^ne^{-\tilde \gamma''(1-\varepsilon)|m-n|^\beta}\leq r^n\}\\ 
 &=&\{n\in\mathbb N \  : \ n\leq \frac{\gamma''(1-\varepsilon)|m-n|^\beta}{\ln(\frac{K}{r})} \}.
 \end{eqnarray*} 
Denote by $n_0$ the highest natural number such that $n_0\in\Lambda_{|m-n|}.$ 
 Then
 
\begin{eqnarray*}
 \sum_{k=1}^{n_0}  |(R^k)_{m,n}|
& \leq& \frac{1}{2P} \sum_{k=1}^{n_0} K^k  e^{-\gamma''(1-\varepsilon)|m-n|^\beta}
 e^{-\varepsilon \gamma''|m-n|^\beta}\\
 &\leq& \frac{1}{2P} e^{-\varepsilon \gamma''|m-n|^\beta}  \sum_{k=1}^\infty  r^k =  \frac{r}{1-r} \frac{1}{2P} e^{-\varepsilon \gamma''|m-n|^\beta}.
\end{eqnarray*}

 For $n>n_0$ we have 
 $r^n<r^{n_0}=e^{n_0\ln r}$
  and hence, 
 $$ \sum_{k=n_0+1}^{\infty}  |(R^n)_{m,n}|
 \leq
 \sum_{k=n_0+1}^{\infty} r^k 
 = r^{n_0+1}\sum_{k=0}^\infty r^k 
< \frac{1}{1-r} e^{-  \gamma''(1-\varepsilon) \frac{\ln (1/r)}{\ln(K/r)} |m-n|^\beta}.$$
 
Let
 $$\gamma_1=
 \min
 \{
 \frac{\ln 1/r}{\ln K/r}{\gamma''(1-\varepsilon), \varepsilon \gamma''}
 \}.
$$
 
Therefore, 
 $$ \sum_{n=0}^\infty  |(R^n)_{m,n}|
\leq  e^{- \gamma_1 |m-n|^\beta  } (1 + \frac{r}{1-r}\cdot \frac{1}{2P} + \frac{1}{1-r}).
$$

Now using the representation  $A^{-1}=A^*(AA^*)^{-1}$ and Lemma \ref{ab}, we can conclude that
\begin{eqnarray*}
|(A^{-1})_{m,n}|
&\leq&  C^A  (1 + \frac{r}{1-r}\cdot \frac{1}{2P} + \frac{1}{1-r}) 
2P_{\gamma-\gamma_1,\beta} 
 e^{-\gamma_1|m-n|^\beta}. 
\end{eqnarray*}

\vspace{.1in}
{\bf Acknowledgements} 
The authors acknowledge support from OeAD GmbH through the  Scientific and Technological  Cooperation projects 
MULT\_DR 01/2017 and SRB 01/2018,
 from the 
Vienna Science and Technology Fund (WWTF) through project VRG12-009, 
and from Project 174024 of the Serbian Ministry of Science. The second author is grateful for the hospitality of the University of Novi Sad, where most of the research on the presented topic was done.

\vspace{.1in}\noindent
Stevan Pilipovi\'c \\
Department of Mathematics and Informatics, University of Novi Sad\\
Trg D. Obradovi\'ca 4, 21000 Novi Sad, Serbia\\
pilipovic@dmi.uns.ac.rs

\vspace{.1in}\noindent
Diana T. Stoeva\\
 Faculty of Mathematics, University of Vienna\\
Oskar-Morgenstern-Platz 1, Vienna 1090, Austria\\
diana.stoeva@univie.ac.at

\end{document}